\documentclass[12pt, reqno]{amsart}

\usepackage{amsmath,amsthm,amsfonts,amssymb}

\usepackage{color}

\theoremstyle{plain}

\newtheorem{thm}{Theorem}[section]
\newtheorem{cor}[thm]{Corollary}
\newtheorem{lem}[thm]{Lemma}
\newtheorem{prop}[thm]{Proposition}

\newtheorem{rem}{Remark}

\numberwithin{equation}{section}

\newcommand{\ann}{\hbox{\rm Ann}}
\definecolor{mjo}{rgb}{0,0,.9}

\newcommand{\cc}{\mathbb{C}}

\newcommand{\z}{\mathbb{Z}}

\newcommand{\nn}{\mathbb{N}}

\newcommand{\setof}[2]{\{ #1 \, | \, #2 \} }

\newcommand{\n}{\mathfrak{n}}

\newcommand{\h}{\mathfrak{h}}

\newcommand{\B}{\mathfrak{b}}

\newcommand{\g}{\mathfrak{g}}

\newcommand{\p}{\mathcal{P}}

\newcommand{\pp}{\tilde{\mathcal{P}}}

\newcommand{\V}{\mathcal{V}}
\def\l{\lambda}

\def\md{\mbox{\rm maxdeg}\,}

\begin{document}

\title{Whittaker Modules for the twisted Heisenberg-Virasoro Algebra}

\author{Dong Liu}
\address{Department of Mathematics, Huzhou Teachers College, Zhejiang Huzhou, 313000, China}
\email{liudong@hutc.zj.cn}
\author{Yuzhu Wu}
\address{Department of Mathematics, Changshu Institute of
Technology, Jiangsu Changshu, 215500, China} \email{yuezhuwujt@126.com}

\author{Linsheng Zhu}
\address{Department of Mathematics, Changshu Institute of
Technology, Jiangsu Changshu, 215500, China}
\email{lszhu@cslg.edu.cn}


\maketitle

\begin{abstract}
We define Whittaker modules for the twisted Heisenberg-Virasoro
algebra and obtain analogues to several results from the classical
setting, including a classification of simple Whittaker modules by
central characters.
\end{abstract}

\section{Introduction}
\setcounter{section}{1}\setcounter{thm}{0}

The twisted Heisenberg-Virasoro algebra has been first studied by
Arbarello et al. in \cite{ADKP}, where a connection is established
between the second cohomology of certain moduli spaces of curves and
the second cohomology of the Lie algebra of differential operators
of order at most one:
$$\mathcal L=\{f(t)\frac{d}{dt}+g(t)|f(t),g(t)\in{\mathbb
C}[t,t^{-1}]\}.$$ Moreover, the twisted Heisenberg-Virasoro algebra
has some relations with the full-toroidal Lie algebras and the $N = 2$
Neveu-Schwarz superalgebra, which is one of the most important
algebraic objects realized in superstring theory. The structure of
the irreducible representations for ${\mathcal L}$  was
studied in \cite{ADKP},\cite{B} and \cite{JJ}. The irreducible weight modules
were classified in \cite{LZ} and \cite{LJ}.

In this paper we study a class of modules for the
twisted Heisenberg-Virasoro algebra.
Whittaker modules were first
discovered for $sl_2(\cc)$ by Arnal and Pinzcon \cite{AP}.  In
\cite{Ko78}, Kostant defined Whittaker modules for an arbitrary
finite-dimensional complex semisimple Lie algebra $\mathfrak{g}$,
and showed that these modules, up to isomorphism, are in bijective
correspondence with ideals of the center $Z(\mathfrak{g})$. In
particular, irreducible Whittaker modules correspond to maximal
ideals of $Z(\mathfrak{g})$. In the quantum setting, Whittaker
modules have been studied for $\mathcal U_h(\mathfrak{g})$ in
\cite{Se00} and for $\mathcal U_q(\mathfrak{sl}_2)$ in \cite{On05}.
Recently Whittaker modules have also been studied by for the
Virasoro algebra in \cite{OW}, for the Schr\"{o}dinger-Virasoro
algebra in \cite{ZT},  for Heisenberg algebras in \cite{Ch08}, for
generalized Weyl algebras in \cite{BO08},  for the $W$-algebra
$W(2,2)$ in \cite{WJ}, and for the super Virasoro
algebra in \cite{LPZ}.

We note that the center of the twisted Heisenberg-Virasoro algebra
is 4-dimensional, which is different from the above algebras whose
centers are all one-dimensional. So our construction is new and it
also offers to construct too much Wittaker modules, although some
calculations of this paper are similar as that in \cite{On05} and
\cite{ZT}.

Throughout the paper, we shall use $\cc, \z,\z_+$ and $\nn$ to denote the sets of the complex numbers, the
integers, the positive integers and the nonnegative integers,
respectively.

\section{Preliminaries} \label{sec:preliminaries}
Let $\V$ denote the twisted Heisenberg-Virasoro Lie algebra. Then
$\V={\rm span}_\cc \{z_1, z_2, z_3, L_k, I_k \mid k \in \z \}$ with
Lie bracket
\begin{align*}
[L_k, L_j] &= (j-k)L_{k+j} + \delta_{j,-k} \frac{k^3-k}{12} z_1; \\
[L_k, I_j] &= jI_{k+j} + \delta_{j,-k}(k^2-k)z_2; \\
[I_k, I_j] &= \delta_{j,-k}k z_3; \\
[z_i, \V] &= 0, i=1,2,3.
\end{align*}
We will make use of the following subalgebras:
\begin{eqnarray*}
\n^+ &= & {\rm span}_\cc \{ L_1, I_1, L_2, I_2, \ldots \} \\
\n^- &= & {\rm span}_\cc \{ L_{-1}, I_{-1}, L_{-2}, I_{-2}, \ldots \} \\
\h&=& {\rm span}_\cc \{ z_1, z_2, z_3, L_0, I_0 \} \\
\B^- &=& \n^- \oplus \h \\
\B^+ &=& \h \oplus \n^+
\end{eqnarray*}
Let $S(z)$ be the symmetric algebra generated by $z_1, z_2,
z_3, z_0=I_0$, that is, polynomials in $z_0, z_1, z_2, z_3$.  Then
$S(z)$ is evidently contained in $Z(\V)$, the center of the
universal enveloping algebra $U(\V)$.
\subsection{Partitions and Pseudopartitions}
The following notation for partitions and pseudopartitions will be used to describe bases for $U(\V)$ and for Whittaker modules.

We define a {\it partition} $\mu$ to be a non-decreasing sequence of positive integers
$\mu=( 0 < \mu_1 \leq \mu_2 \leq \cdots \leq \mu_r)$.  A {\it pseudopartition} $\lambda$ is a non-decreasing sequence of non-negative integers
\begin{equation} \label{eqn:pseudopart1}
\lambda=( 0 \leq \lambda_1 \leq \lambda_2 \leq \cdots \leq \lambda_s).
\end{equation}
Let $\p$ represent the set of partitions, and let $\pp$ denote the set of pseudopartitions.  Then $\p \subseteq \pp$.

We also introduce an alternative notation for partitions and pseudopartitions.  For $\lambda \in \pp$, write
\begin{equation}\label{eqn:pseudopart2}
\lambda = (0^{\lambda(0)}, 1^{\lambda(1)}, 2^{\lambda(2)}, \ldots),
\end{equation}
where $\lambda(k)$ is the number of times $k$ appears in the pseudopartition and $\lambda(k)=0$ for $k$ sufficiently large.   Then a pseudopartition $\lambda$ is a partition exactly when $\lambda(0)=0$.
For $\mu\in \p$, $\lambda \in \pp$, define
\begin{eqnarray*}
| \lambda | &=& \lambda_1 +  \lambda_2 + \cdots + \lambda_s \quad \mbox{(the size of $\lambda$)}\\
\# (\lambda) &=& \lambda(0) + \lambda(1) + \cdots \quad \mbox{(the $\#$ of parts of $\lambda$)}.
\end{eqnarray*}
$$|\mu+\lambda | =|\mu|+|\lambda |$$
For $\mu\in\p, \lambda \in \pp$, define elements $I_{-\mu}, L_{-\lambda} \in U(\V)$ by
\begin{eqnarray*}
I_{-\mu} &=& I_{-\mu_1} I_{-\mu_2}\cdots I_{-\mu_r} = \cdots I_{-2}^{\mu(2)} I_{-1}^{\mu(1)}\\
L_{-\lambda} &=& L_{-\lambda_s} \cdots L_{-\lambda_2}  L_{-\lambda_1}= \cdots L_{-1}^{\lambda(1)} L_0^{\lambda(0)}.
\end{eqnarray*}
Define $\overline{0}= (0^0, 1^0, 2^0, \ldots)$, and write
$L_{\overline{0}}=1 \in U(\V)$. We will consider $\overline{0}$ to
be an element of $\pp$ but not of $\p$. For any $\lambda \in \pp$,
$\mu\in \p$ and $p(z) \in S(z)$, $p(z)L_{- \lambda}I_{-\mu}\in
U(\V)_{-|\lambda+\mu|}$, where $U(\V)_{-|\lambda+\mu|}$ is the
$-|\lambda+\mu|$-weight space of $U(\V)$ under the adjoint action.

\subsection{Whittaker Modules}\label{subsec:WhittakerModules}
\setcounter{equation}{0}
In the classical setting of a finite-dimensional complex semisimple
$\g$, a Whittaker module is defined  in terms of an algebra
homomorphism from the positive nilpotent subalgebra $\g^+$ to $\cc$
(see \cite{Ko78}).  This homomorphism is required to be
\emph{nonsingular}, meaning that it takes nonzero values on the
Chevalley generators of $\g^+$.

In the present setting, the elements $L_1, L_2, I_1\in \n^+$
generate $\n^+$.  Thus we assume that $\psi: \n^+ \rightarrow \cc$
is a Lie algebra homomorphism such that $\psi(L_1), \psi(L_2),
\psi(I_1)\neq 0$. We retain this assumption until to section 5.  The
commutator relations in the definition of $\V$ force
$$
\psi(L_i) = 0\quad\text{for}\  i \ge 3,\quad \psi(I_j) = 0\quad \text{for}\ j\ge 2.
$$
For a $\V$-module $W$, a vector $w \in V$ is a {\it Whittaker
vector} if $xw = \psi(x) w$ for all $x \in \n^+$.  A $\V$-module $W$
is a {\it Whittaker module} if there is a Whittaker vector $w \in W$
which generates $W$.
For a given $\psi: \n^+ \rightarrow \cc$,
define $\cc_\psi$ to be the one-dimensional $\n^+$-module given by
the action $x \alpha = \psi (x) \alpha$ for $x \in \n^+$ and $\alpha
\in \cc$.  The \emph{universal} Whittaker module $M_\psi$ is given
by
$$
M_{\psi}= U(\V) \otimes_{U(\n^+)} \cc_{\psi}.
$$
We use the term universal to refer to the property in Lemma
\ref{lem:universalProp}(i) below.  Let $w=1 \otimes 1 \in M_\psi$.
By the PBW theorem $U( \B^-)$ has a basis $\setof{z^t  L_{-
\lambda}I_{-\mu}}{(\lambda, \mu)\in \pp\times\p, t \in \nn^4}$, where $z^t=z_0^{t_0}z_1^{t_1}z_2^{t_2}z_3^{t_3}$, $t=(t_0,
t_1, t_2, t_3)\in \nn^4$. Thus $M_{\psi}$ has a basis
\begin{equation} \label{eqn:PBWbasis}
\{z^t  L_{- \lambda}I_{-\mu} w \mid (\lambda, \mu)\in \pp\times \p,
t \in \nn^4 \}
\end{equation}
and $uw \neq 0$ whenever $0 \neq u \in U( \B^-)$. Define the {\it
degree}  of $z^t L_{-\lambda}I_{-\mu}w$ to be $| \lambda+\mu |$. For
any $0 \neq v \in M_{\psi}$, define $\md(v)$ to be the maximum
degree of any nonzero component of homogeneous degree, and define
$\md(0) = -\infty$.  Define $\max_{L_0} (v)$ to be the maximum
value $\lambda (0)$ for any term $z^t L_{-\lambda}I_{-\mu}w$ with
nonzero coefficient.

For $\xi=(\xi_0, \xi_1, \xi_2, \xi_3) \in \cc^4$, define
$$L_{\psi, \xi} = M_\psi / \sum_{i=0}^3(z_i - \xi_i1)M_\psi,$$
and let $\overline{\ \cdot \ }: M_{\psi} \rightarrow L_{\psi, \xi}$
be the canonical homomorphism. Similar the Lemma in \cite{OW}, we have

\begin{lem}\label{lem:universalProp}
Fix $\psi$ and $M_\psi$ as above.
\begin{itemize}
\item[(i)] Let  $V$ be a Whittaker module of type $\psi$ with cyclic Whittaker vector $w_V$.  Then there is a surjective map $\varphi: M_{\psi} \rightarrow V$ taking $w=1 \otimes 1$ to $w_V$.
\item[(ii)] Let $M$ be a Whittaker module of type $\psi$ with cyclic Whittaker vector $w_M$, and suppose that for any Whittaker module $V$ of type $\psi$ with cyclic Whittaker vector $w_V$ there exists a surjective homomorphism $\theta: M \rightarrow V$ with $\theta (w_M) = w_V$.  Then $M \cong M_{\psi}$.
\end{itemize}
\end{lem}
\qed

\section{Whittaker vectors in $M_{\psi}$ and $L_{\psi, \xi}$}Fix $\psi : \n^+ \to \cc$ with $\psi(L_1), \psi(L_2),
\psi(I_1) \neq 0$, and let $w = 1 \otimes 1 \in M_\psi$ as before.
In this section we shall determine the Whittaker vectors in the modules $M_{\psi}$
and $L_{\psi, \xi}$.  For
any $w' \in M_{\psi}$, $w'=u w$ for some $u \in U(\B^+)$, we have
$$
(L_n - \psi(L_n)) w' = [L_n, u]w,\quad (I_n - \psi(I_n)) w' = [I_n, u]w.
$$
If $w'$ is a Whittaker vector, then
$$
(L_n - \psi(L_n))w'=(I_n -
\psi(I_n)) w'=0.
$$

\begin{lem}\label{lem:computation2}\cite{On05}
Define $M_\psi$ and $w = 1 \otimes 1 \in M_\psi$ as above, and let $a  \in \z_{+}$, $k  \in \nn$.
Then
$$
[L_{k+2}, L_{-k}^a]w = v-a (2k+2)\psi(L_2) L_{-k}^{a-1} w,
$$
such that $\md(v) <k(a-1)$ if $k>0$, and $\max_{L_0} (v)<a-1$ if
$k=0$.
\end{lem}

\begin{lem}  \ \\
\begin{itemize}
\item[(i)]For $m\in \z_{+},$ we have
$$
\md([I_m,L_{-\lambda}I_{-\mu}]w) \leq|\lambda+\mu|-m+1.
$$

\item[(ii)]$[I_{k+1},L_{-k}^{a}I_{-\mu}]=v-a(k+1)\psi(I_1)I_{-\mu}L_{-k}^{a-1},$
where $\md(v)<(a-1)k$ if $k>0,$ $\max_{L_0}(v)<a-1$ if $k=0.$

\item[(iii)]Suppose
$\lambda=(0^{\lambda(0)},1^{\lambda(1)},2^{\lambda(2)},\cdots)$ and
$k\in \nn$ is the minimal such that $\lambda(k)\neq0,$ then
$$[I_{k+1},L_{-\lambda}I_{-\mu}]w=
v-(k+1)\lambda(k)\psi(I_1)L_{-\lambda{'}}I_{-\mu}w,$$ where if
$k=0$ then $v=v{'}+v{''}$ with $\md(v{'})<|\lambda|-k$ and
 $\max_{L_0}(v{''})<\lambda(k)-1,$ if $k>0$ then
$\md(v)<|\lambda|-k;$ ${\lambda^{'}}$
satisfies $\lambda{'}(k)=\lambda(k)-1, \lambda{'}(i)=\lambda(i)$
for all $i>k.$
\end{itemize}

\end{lem}

\begin{proof}

 (i) Suppose
$\lambda=(\lambda_1,\lambda_2,\cdots,\lambda_r),
\mu=(\mu_1,\mu_2,\cdots,\mu_s).$ Note that $[E_m,-]$ is a derivation
of $U(\V)$ for any $E_m=L_m$ or $I_m$, we have

\begin{eqnarray*}
[I_m,L_{-\lambda}I_{-\mu}]
=&&\sum_{i=1}^{r}L_{-\lambda_r}\cdots[I_m,L_{-\lambda_{i}}]\cdots
L_{-\lambda_1}I_{-\mu}\\
&+&\sum_{j=1}^{s}L_{-\lambda}I_{-\mu_r}\cdots[I_m,I_{-\mu_{j}}]\cdots
I_{-\mu_1}\\
=&&\sum p_{{\lambda{'}},\mu',n}(z_1)L_{-\lambda{'}}I_{-\mu'}I_{n}\\
&+&p'(z_3)L_{-\lambda}I_{-\mu''},
\end{eqnarray*}
where ${\lambda{'}}\in {\mathcal{P}}, n\in \mathbb{Z}$ satisfies
$$-|{\lambda{'}}+\mu'|+n=-|\lambda+\mu|+m=-|\lambda+\mu''|.
$$  Then
$$
\md([I_m,L_{-\lambda}I_{-\mu}]w)=\md(\sum
p_{{\lambda{'}},n}(z)L_{-\lambda{'}}I_{-\mu}I_{n}w)
\leq|\lambda+\mu|-m+1
$$
since $I_nw=\psi(I_n)w=0$ for $n>1.$

(ii)
\begin{eqnarray*}
&&[I_{k+1},L_{-k}^{a}I_{-\mu}]\\
=&&\sum_{i=0}^{a-1}L_{-k}^{i}[I_{k+1},L_{-k}]L_{-k}^{a-i-1}I_{-\mu}+p'(z_3)L_{-\lambda}^aI_{-\mu''}\\
=&&-(k+1)\sum_{i=0}^{a-1}L_{-k}^{i}I_{1}L_{-k}^{a-i-1}I_{-\mu}+p'(z_3)L_{-\lambda}^aI_{-\mu''},
\end{eqnarray*} where $-|\mu|+k+1=-|\mu''|$.

We write $[I_{k+1},L_{-k}^{a}I_{-\mu}]$ a linear combination of
basis of $U(\V)$:
$$[I_{k+1},L_{-k}^{a}I_{-\mu}] =p'(z_3)L_{-\lambda}^aI_{-\mu''}+\sum p_{\lambda,\mu, n}(z)
L_{-\lambda}I_{-\mu}I_{n}-(k+1)ap_{\lambda,\mu,
1}(z)L_{-k}^{a-1}I_{-\mu}I_{1},$$ where if $k=0,$ then $n=1$ and for
each $p_{\lambda,\mu, 1}(z)\neq0$, $L_{-\lambda}=L_{0}^{l}$ for some
$l<a-1;$ if $k>0$ then $n\leq0$ and $\lambda\in{\pp}$ satisfies
$|\lambda|-n=(a-1)k-1.$ Thus (ii) follows.

\vskip 3mm (iii) Denote
$L_{-\lambda}=L_{-\lambda{'}}L_{-k}^{\lambda(k)},$
where $\lambda{'}(i)=\lambda(i)$ for all $i\neq k$ and
$\lambda{'}(k)=0.$ Thus we have
\begin{eqnarray*}
&&[I_{k+1},L_{-\lambda}I_{-\mu}]w
=[I_{k+1},L_{-\lambda{'}}L_{-k}^{\lambda(k)}]I_{-\mu}w+p'(z_3)L_{-\lambda}I_{-\mu''}w\\
&=&[I_{k+1},L_{-\lambda{'}}]L_{-k}^{\lambda(k)}I_{-\mu}w+
L_{-\lambda{'}}[I_{k+1},L_{-k}^{\lambda(k)}]I_{-\mu}w+p'(z_3)L_{-\lambda}I_{-\mu''}w,
\end{eqnarray*} where $-|\mu|+k+1=-|\mu''|$.

We say that the degree of the first and the last summands are
strictly smaller than $|\lambda|-k.$ In fact, if
${\lambda{'}}=\bar{0}$ then
$$
\md([I_{k+1},L_{-\lambda{'}}]L_{-k}^{\lambda(k)}I_{-\mu}w)=-\infty;
$$
If ${\lambda{'}}\neq\bar{0},$ suppose $d$ is the minimal such that
$\lambda{'}(d)\neq0,$ then $d\geq k+1$ by the assumption of $k,$
Therefore
$$
\md([I_{k+1},L_{-\lambda{'}}]L_{-k}^{\lambda(k)}I_{-\mu}w)
=|{\lambda{'}}+\mu|-(k+1)+k\lambda(k)=|\lambda+\mu|-(k+1).
$$ By applying (ii) to the second summand, we complete the proof of
(iii).
\end{proof}

\begin{lem}For any
$(\lambda, \mu)\in \pp\times{\mathcal{P}},$ $m, k\in \mathbb{N},$ we
have
\begin{itemize}
\item[(i)] $\md([L_m,L_{-\lambda}I_{-\mu}]w)
\leq|\lambda+\mu|-m+2;$
\item[(ii)]If $\mu(i)=\lambda(i)=0$ for all
$0\leq i\leq k,$ then
$$\md([L_{k+1},
L_{-\lambda}I_{-\mu}]w) \leq|\lambda+\mu|-k-1;$$
\item[(iii)]If $\mu(j)=\lambda(i)=0$ for all
$0\leq i\leq k,$ $0\leq j<k$ and $\mu(k)\neq0,$ then
$$[L_{k+1},
L_{-\lambda}I_{-\mu}]w=v-k\psi(I_1)\mu(k)L_{-\lambda}I_{-\mu{'}},$$
where $\md(v)<|\lambda+ \mu|-k$ and $\mu{'}$ satisfies that
$\mu{'}(i)=\mu(i)$ for all $i\neq k$ and $\mu{'}(k)=\mu(k)-1.$
\end{itemize}

\end{lem}

\begin{proof}

To prove (i), we write
$[L_m,L_{-\lambda}I_{-\mu}]$ a linear combination of the basis of
$U(\V):$
\begin{eqnarray*}
[L_m,L_{-\lambda}I_{-\mu}]&=& \sum_{\lambda{'},\mu{'}}
p_{\lambda{'},\mu{'}}(z) L_{-\lambda{'}}I_{-\mu{'}}
\\
&+&\sum_{\lambda{''},\mu{''}, n,E_n} p_{\lambda{''},\mu{''},
n}(z)L_{-\lambda{''}}I_{-\mu{''}}E_n,
\end{eqnarray*} where
${\lambda{'}}, \mu{'}, {\lambda{''}}, \mu{''}, n$ satisfy
$-|\mu{'}+{\lambda{'}}|= -|\mu{''}+{\lambda{''}}|+n
=-|\mu+\lambda|+m,$ $E_n=L_n$ or $I_n$. Since $L_iw=I_jw=0$
for $i>2, j>1,$ the first inequality of (i) follows from above
equality. Similarly, we can get the second inequality.

To prove (ii), note first that if $(\lambda,
\mu)=(\bar{0},\underline{0})$ then
$$\md([L_{k+1},
L_{-\lambda}I_{-\mu}]w)=-\infty,$$ we suppose $(\lambda,
\mu)\neq(\bar{0},\underline{0})$ and write
$$L_{-\lambda}I_{-\mu}
=L_{-\lambda_r}\cdots L_{-\lambda_1}I_{-\mu_s}\cdots I_{-\mu_1}.$$
Since $k<\mu_j, \lambda_i$ for all $1\leq j\leq s, 1\leq a\leq
t,1\leq i\leq r,$ we have
$$[L_{k+1},
L_{-\lambda}I_{-\mu}]\in U(\mathfrak{b}^{-}),$$ and
$$\md([L_{k+1},
L_{-\lambda}I_{-\mu}]w) =|\mu+\lambda|-k-1.$$ So (ii) holds.

With respect to (iii), noting
\begin{eqnarray*}[L_{k+1},
L_{-\lambda}I_{-\mu}]
&=&L_{-\lambda}[L_{k+1},I_{-\mu^{'}}]I_{-k}^{\mu(k)}
\\
&&+L_{-\lambda}I_{-\mu^{'}}[L_{k+1},I_{-k}^{\mu(k)}]
\\
&&+[L_{k+1}, L_{-\lambda}]I_{-\mu^{'}}I_{-k}^{\mu(k)},
\end{eqnarray*}
where $\mu{'}\in \mathcal{P}$ such that $\mu{'}(k)=0,
\mu{'}(i)=\mu(i)$ for $i\neq k.$ By assumption about $k,$ we see
that $[L_{k+1},I_{-\mu{'}}],$ $[L_{k+1},
L_{-\lambda}]\in
U(\mathfrak{b}^{-}),$
$[L_{k+1},I_{-k}^{\mu(k)}]=-k\mu(k)I_{-k}^{\mu(k)-1}I_1,$ and
$$\md(L_{-\lambda}[L_{k+1},I_{-\mu{'}}]I_{-k}^{\mu(k)}w)
\leq |\mu+\lambda|-k-1,$$
$$
\md([L_{k+1},
L_{-\lambda}]I_{-\mu{'}}I_{-k}^{\mu(k)}w)\leq | \mu+\lambda|-k-1,$$
$$L_{-\lambda}I_{-\mu{'}}[L_{k+1},I_{-k}^{\mu(k)}]
w=-k\mu(k)\psi(I_1) L_{-\lambda}I_{-\mu{''}}w,$$ where $\mu{''}\in
\mathcal{P}$ such that $\mu{''}(k)=\mu(k)-1, \mu{''}(i)=\mu(i)$ for all
$i\neq k.$ Thus (iii) holds.
\end{proof}

\begin{prop}\label{prop:describeWhitVectsM}
Let $M_{\psi}$ be a universal Whittaker module for $\V$, generated by the Whittaker vector $w= 1 \otimes 1 $.  If $w' \in M_{\psi}$ is a Whittaker vector, then $w'=p(z)w$ for some $p(z) \in S(z)$.
\end{prop}
\begin{proof}
Let $w' \in M_{\psi}$ be an arbitrary vector.  By
(\ref{eqn:PBWbasis}),
$$w' = \sum_{(\lambda, \mu) \in (\pp, \p)} p_{\lambda, \mu}(z)  L_{-\lambda}I_{-\mu}w$$
 for some polynomials $p_{\lambda, \mu}(z) \in S(z)$.
 We will show that  if there is $\lambda\neq \overline{0}=(0^0, 1^0, 2^0, \ldots )$ or $\mu\neq \underline{0}=(1^0, 2^0, \ldots )$ such that
 $p_{\lambda, \mu}(z) \neq 0$, then there is $m \in \z_{>0}$ such that $(L_m - \psi(L_m))w'  \neq 0$
 or $(I_m - \psi(I_m))w'  \neq 0$.  In this case $w'$ is not a Whittaker vector, which proves the result.

Let $N= \mbox{max} \{ | \lambda+\mu | \mid  p_{\lambda, \mu}(z) \neq 0 \}$,
and define $\Lambda_N = \{ (\lambda, \mu) \in \pp\times\p \mid  p_{\lambda, \mu}(z) \neq 0, |\lambda+\mu| = N \}$.

Assume $N>0.$ Set
$$k:=min\{n\in \mathbb{Z}_{+}|
\mu(n)\neq0 \ or \ \lambda(n)\neq0 \ for \  some \ (\lambda, \mu)\in
\Lambda_N\}.$$

\vskip 3mm\noindent{\bf{Case I.}} $k$ satisfies $\lambda(k)\neq0$
for some $(\lambda, \mu)\in \Lambda_N.$ By Lemma 3.1 , we have
\begin{eqnarray*}
&&(I_{k+1}-\psi(I_{k+1}))w{'}\\
&=&\sum_{(\lambda, \mu)\notin\Lambda_N}
p_{\lambda, \mu}(z)[I_{k+1},
L_{-\lambda}I_{-\mu}]w\\
&&+\sum_{\begin{subarray}\ (\lambda, \mu)\in\Lambda_N\\
\ \ \lambda(k)=0
\end{subarray}} p_{\lambda, \mu}
(z)[I_{k+1},
L_{-\lambda}I_{-\mu}]w\\
&&+\sum_{\begin{subarray}\ (\lambda, \mu)\in\Lambda_N\\
\ \ \lambda(k)\neq0
\end{subarray}} p_{\lambda, \mu}
(z)[I_{k+1}, L_{-\lambda}I_{-\mu}]w
\end{eqnarray*}
By using Lemma 3.2 (i) to the first summand, we know that the
degree of it is strictly smaller than $N-k.$ As for the second
summand, note that $\lambda(i)=0$ for $0\leq i\leq k,$ we have
$$[I_{k+1},L_{-\lambda}I_{-\mu}]
=[I_{k+1},L_{-\lambda}]I_{-\mu}+L_{-\lambda}[I_{k+1},I_{-\mu}]
\in\mathcal{S}(z)U(\V^{-}),$$ thus the degree of it is also strictly
smaller than $N-k.$ Now using Lemma 3.2 (iii) to the third summand,
we know it has form:
$$v-\sum_{\begin{subarray}
\ \ (\lambda, \mu)\in\Lambda_N\\
\ \ \lambda(k)\neq0
\end{subarray}}(k+1)\lambda(k)\psi(I_1)p_{\lambda, \mu}
(z)L_{-\lambda}I_{-\mu}w,$$ where if $k=0$ then $v=v{'}+v{''}$
such that $\md(v{'})<N-k$ and $\max_{L_0}(v{''})<\lambda(0)-1;$
if $k>0$ then $\md(v)<N-k.$ ${\lambda{'}}$ satisfies
$\lambda{'}(k)=\lambda(k)-1,\lambda{'}(i)=\lambda(i)$ for all
$i>k.$ Thus the degree of the third summand is $N-k$, which implies
$(I_{k+1}-\psi(I_{k+1}))w{'}\neq0.$

\vskip 3mm\noindent{\bf{Case II.}} $k$ satisfies $\mu(k)\neq0,
\lambda(k)=0$ for some $(\lambda, \mu)\in \Lambda_N.$

For convenience, we denote $w{'}=w{'}_1+w{'}_2$ such that
$$
w{'}_1=\sum_{(\lambda, \mu)\notin\Lambda_N}
p_{\lambda, \mu} (z) L_{-\lambda}I_{-\mu}w,
$$
$$
w{'}_2=\sum_{(\lambda, \mu)\in\Lambda_N}
p_{\lambda, \mu} (z) L_{-\lambda}I_{-\mu}w.
$$ We will find some
element $E_m\in \V^+$ such that $(E_m-\psi(E_m))w{'}\neq0$
according to the following subcases:

\vskip 3mm\noindent{\bf{Subcase 1.}} $\md(w{'}_1)<N-1.$ In this
subcase,
\begin{eqnarray*}
(L_{k+1}-\psi(L_{k+1}))w'_{1} =\sum_{(\lambda,
\mu)\notin\Lambda_N} p_{\lambda, \mu} (z)[L_{k+1},
L_{-\lambda}I_{-\mu}]w.
\end{eqnarray*}
By the first inequality of Lemma 3.3 (i), we have
$$
\md((L_{k+1}-\psi(L_{k+1}))w'_{1})<(N-1)-(k+1)+2=N-k.
$$
\begin{eqnarray*}
(L_{k+1}-\psi(L_{k+1}))w'_{2} & = &\sum_{\begin{subarray}
\ (\lambda, \mu)\in\Lambda_N\\
\ \ \ \mu(k)=0\end{subarray}}
p_{\lambda, \mu} (z)[L_{k+1},
L_{-\lambda}I_{-\mu}]w\\
&&+\sum_{\begin{subarray}
\ (\lambda, \mu)\in\Lambda_N\\
\ \ \ \mu(k)\neq0\end{subarray}} p_{\lambda, \mu} (z)[L_{k+1},
L_{-\lambda}I_{-\mu}]w.
\end{eqnarray*}
By using Lemma 3.3 (ii) to the first summand  and (iii) to the
second, we have
$$(L_{k+1}-\psi(L_{k+1}))w^{'}_{2}=v-\sum_{\begin{subarray}
\ (\lambda, \mu)\in\Lambda_N\\
\ \ \ \mu(k)\neq0\end{subarray}}
p_{\lambda, \mu} (z)\mu(k)k\psi(I_1)
I_{-\mu{'}}L_{-\lambda}w,$$
where $\md(v)<N-k,$ $\mu{'}$ satisfies that $\mu{'}(i)=\mu(i)$
for all $i\neq k$ and $\mu{'}(k)=\mu(k)-1.$ Thus, we have
$$(L_{k+1}-\psi(L_{k+1}))w{'}=
(L_{k+1}-\psi(L_{k+1}))w'_{1}+(L_{k+1}-\psi(L_{k+1}))w'_{2}\neq0$$

\vskip 3mm\noindent{\bf{Subcase 2.}} $\md(w{'}_1)=N-1$ and
$\lambda=\bar{0}$ for any
$(\lambda, \mu)$ with
$p_{\lambda, \mu}(z)\neq0.$ In this
subcase,
$$w'_1=\sum_{(\bar{0}, \mu)\notin\Lambda_N}
p_{\bar{0}, \mu} (z)
I_{-\mu}w,$$
$$w'_2=\sum_{(\bar{0}, \mu)\in\Lambda_N}
p_{\bar{0}, \mu}(z)
I_{-\mu}w.$$
Hence
\begin{eqnarray*}
(L_{k+1}-\psi(L_{k+1}))w'_{1}
=\sum_{(\bar{0}, \mu)\notin\Lambda_N}
p_{\bar{0}, \mu}(z)[L_{k+1},
I_{-\mu}]w,
\end{eqnarray*}
and by the second inequality of Lemma 3.3 (i), we have
$$
\md((L_{k+1}-\psi(L_{k+1}))w'_{1})<N-k;
$$
\begin{eqnarray*}
(L_{k+1}-\psi(L_{k+1}))w^{'}_{2} &=&\sum_{\begin{subarray} \
(\bar{0}, \mu)\in\Lambda_N\\
\ \ \ \mu(k)=0
\end{subarray}}
p_{\bar{0}, \mu}(z)[L_{k+1},
I_{-\mu}]w\\
&+& \sum_{\begin{subarray} \
(\bar{0}, \mu)\in\Lambda_N\\
\ \ \ \mu(k)\neq0
\end{subarray}}
p_{\bar{0}, \mu}(z)[L_{k+1},
I_{-\mu}]w
\end{eqnarray*}
and by Lemma 3.3 (ii) and (iii),
$$
(L_{k+1}-\psi(L_{k+1}))w'_{2}=v-\sum_{\begin{subarray} \
(\bar{0}, \mu)\in\Lambda_N\\
\ \ \ \mu(k)\neq0
\end{subarray}}
p_{\bar{0}, \mu}(z)k\mu(k)
\psi(I_1)I_{-\mu{'}}w,$$ where
$\md(v)<N-k$ and $\mu{'}$ satisfies $\mu{'}(i)=\mu(i)$ for all
$i\neq k$ and $\mu{'}(k)=\mu(k)-1.$ Thus
$\md((L_{k+1}-\psi(L_{k+1}))w)=N-k$ and
$(L_{k+1}-\psi(L_{k+1}))w\neq0.$

\vskip 3mm\noindent{\bf{Subcase 3.}} $\md(w'_1)=N-1,$
$\lambda=\bar{0}$ for any $(\lambda, \mu)\in\Lambda_N$ and
$\lambda\neq\bar{0}$ for some $(\lambda, \mu)$ such that $(\lambda,
\mu)\notin\Lambda_N$ with $p_{\lambda, \mu}(z)\neq0.$ Suppose that
$l$ is the minimal such that $\lambda(l)\neq0$ for some $(\lambda,
\mu)$ which satisfies $(\lambda, \mu)\notin\Lambda_N$ and
$p_{\lambda, \mu}(z)\neq0.$ Let $N^{'}=\{max|\mu+\lambda| |(\lambda,
\mu)\notin\Lambda_N, \lambda(l)\neq0, p_{\lambda, \mu}(z)\neq0\}.$
Then by Lemma 3.2 (iii) we have
$$(I_{l+1}-\psi(I_{l+1}))w'_1
=v-\sum_{\begin{subarray}
\ \ \ \ \  (\lambda, \mu)\notin\Lambda_N\\
\ |\mu+|\lambda|=N{'}\\
\ \ \ \ \ \  \lambda(l)\neq0\end{subarray}}
(l+1)\lambda(l)\psi(I_{1})p_{\lambda, \mu}(z)
I_{-\mu}L_{{\lambda{'}}}w,$$
where $\md(v)<N'-l,$ $ \lambda{'}$ satisfies that
$\lambda{'}(i)=\lambda(i)$ for all $i\neq l$ and
$\lambda{'}(l)=\lambda(l)-1.$ It is clear that $$(I_{l+1}-\psi(I_{l+1}))w'_2=0$$ since
$w'_2=\sum_{(\bar{0}, \mu)\in\Lambda_N}
p_{\bar{0}, \mu}(z)
I_{-\mu}w$ and
$[I_{l+1},I_{-\mu}]=0.$  Thus
$$
(I_{l+1}-\psi(I_{l+1}))w{'}\neq0.
$$

\vskip 3mm\noindent{\bf{Subcase 4.}} $\md(w'_1)=N-1$ and there
exists some $(\lambda, \mu)\in\Lambda_N$
such that $\lambda\neq\bar{0}.$ In this subcase we
suppose $l$ is the minimal such that $\lambda(l)\neq0$ for some
$(\lambda, \mu)\in\Lambda_{N}.$ Then
$l>k$ according to the assumption about $k$. Since
\begin{eqnarray*}
(I_{l+1}-\psi(I_{l+1}))w^{'}_{1} =\sum_{(\lambda,
\mu)\notin\Lambda_N} p_{\lambda, \mu} (z)[I_{l+1},
L_{-\lambda}I_{-\mu}]w,
\end{eqnarray*}
by Lemma 3.2 (i) we have
$$
\md((I_{l+1}-\psi(I_{l+1}))w'_{1})\leq(N-1)-l.
$$
By Lemma 3.2 (iii), we have
\begin{eqnarray*}
(I_{l+1}-\psi(I_{l+1}))w'_{2}
&=&\sum_{(\lambda, \mu)\in\Lambda_N}
p_{\lambda, \mu} (z)[I_{l+1},
L_{-\lambda}I_{-\mu}]w\\
&=&v-\sum_{\begin{subarray}
\ \ \   (\lambda, \mu)\in\Lambda_N\\
\ \ \ \   \lambda(l)\neq0\end{subarray}}
p_{\lambda, \mu}
(z)(l+1)\psi(I_1)\lambda(l)
I_{-\mu}L_{-\lambda'}w,
\end{eqnarray*}
where $\md(v)<N-l,$ $ \lambda'$ satisfies that
$\lambda'(i)=\lambda(i)$ for all $i\neq l$ and
$\lambda'(l)=\lambda(l)-1.$ Thus, we have
$$
\md((I_{l+1}-\psi(I_{l+1}))w'_2)=N-l$$
and
$$
(I_{l+1}-\psi(I_{l+1}))w'=
(I_{l+1}-\psi(I_{l+1}))w'_{1}+(I_{l+1}-\psi(I_{l+1}))w'_{2}\neq0.
$$

\end{proof}

\begin{prop}\label{Prop:SimpleWhittaker}
Let $w = 1 \otimes 1 \in M_{\psi}$ and $\overline{w}= \overline{1 \otimes 1} \in L_{\psi, \xi}$.  If $w' \in L_{\psi, \xi}$ is a Whittaker vector, then $w'=c\overline{w}$ for some $c \in \cc$.
\end{prop}

\begin{proof}
Note that the set $\{L_{- \lambda}I_{-\mu} \overline{w} \mid (
\lambda, \mu) \in (\pp, \p) \}$ is a basis for $L_{\psi, \xi}$. With
this fact now established, it is possible to use essentially the
same argument as in Proposition \ref{prop:describeWhitVectsM} to
complete the proof.  However, we simply replace the polynomials
$p_\lambda (z)$ in $z$ with scalars $p_\lambda$ whenever necessary.
\end{proof}

\section{Simple Whittaker Modules}
By Proposition \ref{Prop:SimpleWhittaker}, we shall  show that
the modules $L_{\psi,\xi}$ are simple and form a complete set of
simple Whittaker modules, up to isomorphism.

Fix an algebra homomorphism $\psi: \n^+ \rightarrow \cc$ and a
Whittaker module $V$ of type $\psi$. We may regard $V$ as an
$\n^+$-module by restriction. Define a modified action of $\n^+$ on
$V$ (denoted by $\cdot$) by setting $x \cdot v = xv - \psi (x) v$
for $x \in \n^+$ and $v \in V$.  Thus if we regard a Whittaker
module $V$ as an $\n^+$-module under the dot action, it follows that
$E_n \cdot v = E_nv - \psi(E_n)v$ for $E_n=L_n$ or $I_n$, $n>0$ and
$v \in V$.

\begin{lem}\label{lem:locallyNilpotent}
  If $n>0$, then $L_n$ and $I_n$ are locally nilpotent on $V$ under the dot action.
\end{lem}
\begin{proof}
By direction calculation.
\end{proof}

By straightfordwork calculation, we have the following lemma
\begin{lem}\label{lem:predotActionFinite}
Let $\lambda \in \pp$ and $t\in\nn^4$.
\begin{itemize}
\item[(i)] For all $n>0$, $E_n \cdot (z^t L_{-\lambda}I_{-\mu}w) \in
\mbox{span}_{\cc} \{ z^{t'} L_{-\lambda'}I_{-\mu'} w \mid |
\mu'+\lambda'| + \lambda'(0) \leq |\mu+\lambda| + \lambda(0)
  \}$.
\item[(ii)]  If $n > |\lambda+\mu|+2$, then $E_n \cdot (L_{-\lambda}I_{-\mu}w) = 0$.
\end{itemize}
\end{lem}\qed

\begin{cor}\label{lem:dotActionFinite}
Suppose $V$ is a Whittaker module for $\V$, and let $v \in V$.
Regarding $V$ as an $\n^+$-module under the dot action, $U(\n^+)
\cdot v$ is a finite-dimensional submodule of $V$.
\end{cor}
\begin{proof}
This is the direct result of Lemma \ref{lem:predotActionFinite}.
\end{proof}

\begin{thm} \label{thm:submodule}
Let $V$ be a Whittaker module for $\V$, and let $S \subseteq V$ be a submodule.  Then there is a nonzero Whittaker vector $w' \in S$.
\end{thm}
\begin{proof}
Regard $V$ as an $\n^+$-module under the dot-action. Let $0 \neq v
\in S$, and let $F$ be the submodule of $S$ generated by $v$ under
the dot-action of $\n^+$.  By Lemma  \ref{lem:dotActionFinite}, $F$
is a finite-dimensional $\n^+$-module.  Since Lemma
\ref{lem:predotActionFinite} implies that $E_n \cdot F = 0$ for
sufficiently large $n$, the quotient of $\n^+$ by the kernel of this
action is also finite-dimensional.  Note that $E_n$ is locally
nilpotent on $V$ (and thus on $F$) under this action by Lemma
\ref{lem:locallyNilpotent}.  Thus Engel's Theorem implies that there
exists a nonzero $w' \in F \subseteq S$ such that $x \cdot w'=0$ for
all $x \in \n^+$.  By definition of the dot-action, $w'$ is a
Whittaker vector.
\end{proof}

\begin{cor} \label{cor:L_PsiXiSimple}
For any $\xi \in \cc^4$, $L_{\psi, \xi}$ is simple.
\end{cor}
\begin{proof}
Let $S$ be a nonzero submodule of $L_{\psi, \xi}$. Since $z \in \V$
acts by the scalar $\xi$ on $L_{\psi, \xi}$, it follows from Theorem
\ref{thm:submodule} that there is a nonzero Whittaker vector $w' \in
S$.  Proposition \ref{Prop:SimpleWhittaker} implies that $w' = cw$
for some $c \in \cc$, and therefore $w \in S$.  Since $w$ generates
$L_{\psi, \xi}$, we have $S= L_{\psi, \xi}$.
\end{proof}

\noindent
\begin{cor} \label{cor:SimpleL_PsiXi}
Let $\psi : \n^+ \to \cc$ be a Lie algebra homomorphism and
$\psi (L_1), \psi (L_2), \psi(I_1) \neq 0$. Let $S$ be a simple
Whittaker module of type $\psi$ for $\V$. Then $S \cong L_{\psi,
\xi}$ for some $\xi \in \cc^4$.
\end{cor}
\begin{proof}
Let $w_s \in S$ be a cyclic Whittaker vector corresponding to
$\psi$.  By Schur's lemma, the center of $U( \V )$ acts by a scalar,
and this implies that there exists $\xi \in \cc^4$ such that $z_i s
= \xi_i  s, i=0,1,2,3$ for all $s \in S$.  Now by the universal
property of $M_\psi$, there exists a module homomorphism $\varphi :
M_\psi \to S$ with $uw \mapsto uw_s$.  This map is surjective since
$w_s$ generates $S$.  But then
$$\varphi \left( \sum_{i=0}^3(z_i -
\xi_i1) M_\psi \right) = \sum_{i=0}^3(z_i - \xi_i1) \varphi (M_\psi
) = \sum_{i=0}^3(z_i - \xi_i1)S = 0.$$
It follows that
$$
\sum_{i=0}^3(z_i - \xi_i1) M_\psi \subseteq \ker \varphi \subseteq M_\psi.
$$
Because $L_{\psi, \xi}$ is simple and $\ker \varphi \neq M_\psi$,
this implies $\ker \varphi = \sum_{i=0}^3(z_i - \xi_i1) M_\psi$.
\end{proof}

For a given $\psi : \n^+ \to \cc$ and $\xi \in \cc^4$, note that $L
= \sum_{i=0}^3 U(\V)(z_i - \xi_i1) + \sum_{m>0} U(\V)(L_m -
\psi(L_m)1) +\sum_{m>0} U(\V)(I_m - \psi(I_m)1)\subseteq U(\V)$ is a
left ideal of $U(\V)$. For $u \in U(\V)$, let $\overline u$ denote
the coset $u + L \in U(\V)/L$. Then we may regard $U(\V)/L$ as a
Whittaker module of type $\psi$ with cyclic Whittaker vector
$\overline 1$.

\begin{lem}\label{lem:quotientEnveloping}
Fix $\psi : \n^+ \to \cc$ with $\psi(L_1), \psi(L_2), \psi(I_1) \neq
0$. Define the left ideal $L$ of $U(\V)$ by $L = \sum U(\V)(z_i -
\xi_i1) + \sum_{m>0} U(\V)(L_m - \psi(L_m)1)+\sum_{m>0} U(\V)(I_m -
\psi(I_m)1)$, and regard $V = U(\V) / L$ as a left $U(\V)$-module.
Then $V$ is simple, and thus $V \cong L_{\psi, \xi}$.
\end{lem}
\begin{proof}
Note that the center of $U( \V )$ acts by the scalar $\xi$ on $V$.
By the universal property of $M_\psi$, there exists a module
homomorphism $\varphi : M_\psi \to V$ with $uw \mapsto u\overline
1$.  This map is surjective since $\overline 1$ generates $V$.  But
then $\varphi \left( \sum(z_i - \xi_i1) M_\psi \right) = \sum(z_i -
\xi_i1) \varphi (M_\psi ) = \sum(z_i - \xi_i1)V = 0$, so it follows
that
$$\sum(z_i - \xi_i1) M_\psi \subseteq \ker \varphi \subseteq M_\psi.$$
Because $L_{\psi, \xi}$ is simple and $\ker \varphi \neq M_\psi$,
this forces $\ker \varphi = \sum(z_i - \xi_i1) M_\psi$.
\end{proof}

\begin{prop}\label{Prop:zScalarImpliesSimple}
Suppose that $V$ is a Whittaker module of type $\psi$ such that $z_i
\in \V$ acts by the scalar $\xi_i \in \cc, i=0,1,2,3$.   Then $V$ is
simple. Moreover, if $w$ is a cyclic Whittaker vector for $V$,
$\ann_{U(\V)} (w) = \sum_{i=0}^3U(\V) (z_i-\xi_i1) + \sum_{m>0}
U(\V)(L_m - \psi(L_m)1)+\sum_{m>0} U(\V)(I_m - \psi(I_m)1)$.
\end{prop}
\begin{proof}
Let $K$ denote the kernel of the natural surjective map $U(\V) \to
V$ given by $u \mapsto uw$.  Then $K$ is a proper left ideal
containing $L =\sum_{i=0}^3U(\V) (z_i-\xi_i1)  + \sum_{m>0}
U(\V)(L_m - \psi_m1)+\sum_{m>0} U(\V)(I_m - \psi(I_m)1)$. By Lemma
\ref{lem:quotientEnveloping}, $L$ is maximal, and thus $K=L$ and $V
\cong U(\V)/L$ is simple.
\end{proof}

\begin{rem}
By Schur's Lemma, Proposition \ref{Prop:zScalarImpliesSimple} applies to any simple Whittaker module.

\end{rem}

\section{Whittaker modules for the differential operator algebra of degree no more than 1} \label{sec:witt}

In this section, we describe the Whittaker modules for the Lie
subalgebra $\mathcal L$ of the differential operator algebra of
degree no more than 1.  Recall that the twisted Heisenberg-Virasoro
algebra $\V$ is the universal central extension of $\mathcal L$.  We
abuse notation and regard $\mathcal L={\rm span}_\cc \{ L_k, I_k
\mid k \in \z \}$ with Lie bracket given by
\begin{align*}
[L_k, L_j] &= (j-k)L_{k+j}, \ [I_k, I_j] = 0,\ [L_k, I_j] = jI_{k+j}
\end{align*}
for $j, k \in \z$.

As $\V$ is the universal central extension of $\mathcal L$,  there
is a surjective Lie algebra homomorphism $\rho : \V \to \mathcal L$
with $\ker \rho=\cc\{z_1,z_2,z_3\}$.  This map extends to a
surjective homomorphism $U(\V) \to U(\mathcal L)$ which we also
denote by $\rho$.

Define the subalgebra $\n^+_\mathcal L \subseteq \mathcal L$ in the
obvious manner.   Since $\n^+ \cong \n^+_\mathcal L$, we make no
distinction between a homomorphism $\psi : \n^+_\mathcal L \to \cc$
and a homomorphism $\psi : \n^+ \to \cc$.  Let $\psi: \n^+_\mathcal
L \rightarrow \cc$ be a Lie algebra homomorphism such that
$\psi(L_1), \psi(L_2), \psi(I_1)\neq 0$.  A $\mathcal L$-module $V$
is a {\it Whittaker module} if there is some $w \in V$ such that $w$
generates $V$ and $xw=\psi(x)w$ for all $x \in \n^+_\mathcal L$.

\begin{prop}
Fix a homomorphism $\psi : \n^+_\mathcal L \to \cc$ with $\psi
(L_1)$, $\psi (L_2)$, $\psi(I_1) \neq 0$,  and let $V$ be a nonzero
Whittaker module of type $\psi$ for $\mathcal L$.  Then $V$ is
simple. Moreover, $V \cong L_{\psi,\eta_0}$ when $L_{\psi, {\bf
0}}$ is viewed as a $\mathcal L$-module, where $\eta_0=(\xi_0, 0,
0, 0)$ for some $\xi_0\in\cc$.
\end{prop}
\begin{proof}
Let $V_\V$ be the $\V$-module obtained by letting $x \in \V$ act on
$V$ by $\rho (x) \in \mathcal L$.  Then $V_\V$ is a nonzero
Whittaker module for $\V$, and the central element $z_1, z_2, z_3
\in \V$ act by $0$s. By Proposition \ref{Prop:zScalarImpliesSimple},
we then have $V_{\V} \cong L_{\psi,\eta_0}$ for some $\eta_0=(\xi_0, 0, 0, 0)\in\cc^4$. As $V_\V$ is the pullback of $V$ and
$\rho : \V \to \mathcal L$ is surjective, we conclude that $V$ must
also be simple.

To check that $L_{\psi, \eta_0}$ can be viewed at a $\mathcal
L$-module, we note that $z_1, z_2, z_3 \in \ann_{\V} (L_{\psi,
\eta_0})$. Therefore, the action of $\mathcal L = \V/
\cc\{z_1,z_2, z_3\}$ is well-defined. Since $V_{\V} \cong
L_{\psi,\eta_0}$ as $\V$ modules, we must have $V \cong
L_{\psi,\eta_0}$ as $\mathcal L$-modules.
\end{proof}

\section{ A class of new module similar to
Whittaker modules}

In this section we construct and study a family of $\V$-module
$V_{\psi,\xi}$ which are generated by $\psi$-vectors. When $\psi$ is
identically zero, each of these modules has a quotient module
isomorphic to a Verma module. When $\psi$ is non-singular these are
the irreducible Whittaker modules $L_{\psi,\xi}$ desired above.
Otherwise $V_{\psi,\xi}$ will be new.

Let $\xi=(\xi_0,\xi_1,\xi_2,\xi_3)\in\cc^4$ and let
$\psi:\V^+\rightarrow \cc$ be any Lie algebra homomorphism. View
$\cc_{\psi,\xi}=\cc w$ as a one-dimensional ${\mathcal Z}+\V^+$-module by
$$xw=\psi(x)w,\,z_0w=\xi_1w,\,z_1w=\xi_2w,\,z_2w=\xi_3w,\,z_3w=\xi_4w,$$
where $x\in\V^+,\,{\mathcal Z}=span_\cc\{z_i|i=0,1,2,3\}$. Then we have an induce
$U(\V)$-module
$$V_{\psi,\xi}:=U(\V)\otimes_{U(\V^+)\mathcal{S}(z)}C_{\psi,\xi}.$$
Clearly, $V_{\psi,\xi}$ has the following facts:

(1) If $\psi$ is non-singular, $V_{\psi,\xi}$ is the irreducible
Whittaker module $L_{\psi,\xi}$.

(2) If $\psi=0$ and $V$ is the submodule generared by
$(L_0-\zeta)w,\,\zeta\in\cc
,$ then the quotient module
$V_{\psi,\xi}/V$ is the Verma module for $\V$.

 From now on we fix a Lie algebra homomorphism $\psi$ such that
 $\psi$ is singular but not identically zero.
 \begin{thm}
For $\xi=(\xi_1,\xi_2,0,0)\in\cc^4$, $V_{\psi,\xi}$ is simple if and
only if $\psi(I_1)\ne 0$.
 \end{thm}
\begin{proof} If $\psi(I_1)=0$, the submodule $V$ of $V_{\psi,0}$
generated by $\{I_{-n}v|n\in\z_+\}$ is proper since $w\notin V$.

Now we suppose $\psi(I_1)\ne 0$. For any $0\ne v\in V$, Write
$$v=\sum a_{\mu,\l,i}(z)L_{\l}I_{-\mu}w,$$
where $a_{\mu,\l,i}(z)$ is in the polynomial algebra of $z_1,\,z_2$.
Set
\begin{eqnarray*}
 N&=& \mbox{max} \{ | \lambda+\mu | \mid  a_{
\lambda,\mu,i}(z) \neq 0 \},\\
\Lambda_N &=& \{ (\lambda,\mu ) \in \pp\times\p \mid a_{\lambda,
\mu,i}(z) \neq 0, |\lambda+\mu| = N \}.
\end{eqnarray*}Now we will
use induction on $N$ to prove that $w\in V.$ If $N=0$, then
$$v=\sum a_{i}(z)L_0^iw.$$ By applying the dot action of $I_1$ on $v$ and using induction on $i$, one can
 prove $w\in V$.

Now assume $N>0$.  Using the same argument of Proposition
\ref{prop:describeWhitVectsM}, one can prove that there is $m \in
\z_{>0}$ such that $L_m \cdot v  \neq 0$, $\,{\rm maxdeg} (L_m \cdot
v)<N$
 or $I_m\cdot v  \neq 0$,$\,{\rm maxdeg} (I_m \cdot v)<N$.
\end{proof}
\vskip30pt \centerline{\bf ACKNOWLEDGMENTS}

\vskip15pt Project is supported by the NNSF (Grant 10671027,
10701019, 10826094), the ZJZSF(Grant Y607136), Qianjiang
Excellence Project of Zhejiang Province (No. 2007R10031), and NSF
08KJD110001 of Jiangsu Educational Committee.

\end{document}